\documentclass[11pt,epsfig]{article}
\usepackage{amssymb}
\usepackage{amsmath}
\usepackage{amsthm}
\usepackage{epsfig}
\usepackage{hyperref}
\usepackage{amsfonts}
\usepackage{color}

\newtheoremstyle{thmm}{1.5ex plus 1ex minus .2ex}{1.5ex plus 1ex minus
.2ex}{\rmfamily}{}{\bfseries}{}{1em}{}
\theoremstyle{thmm}
\newtheorem{theorem}{Theorem}[section]
\newtheorem{lemma}{Lemma}[section]

\renewenvironment{proof}[1][Proof]{\noindent\textit{#1. } }{\hfill$\square$=
}

\allowdisplaybreaks \textwidth  6.0in \textheight 9in
\newcommand{\nn}{\nonumber}
\def \endproof{\vrule height8pt width 5pt depth 0pt}
\def\refe#1{(\ref{#1})}

\def\d{\delta}

\def\R{\mathbb{R}}

\def\d{\,{\rm d}}

\begin{document}

\date{}

\title{\bf Error analysis of linearized semi-implicit Galerkin
finite element methods for nonlinear parabolic equations}

\author{Buyang
Li\,\setcounter{footnote}{0}\footnote{Department
of Mathematics, Nanjing University, Nanjing, P.R. China.
{\tt buyangli@nju.edu.cn} }~\footnotemark[2]
~~and ~Weiwei Sun\footnote{Department
of Mathematics, City University of Hong Kong,
Kowloon, Hong Kong.
The work of the
authors was supported in part by a grant
from the Research Grants
Council of the Hong Kong Special Administrative
Region, China
(Project No. CityU 102005)  {\tt
maweiw@math.cityu.edu.hk} }}

\maketitle

\begin{abstract}
This paper is concerned with the time-step condition
of commonly-used linearized semi-implicit schemes for nonlinear
parabolic PDEs with Galerkin finite element approximations. In
particular, we study the time-dependent nonlinear Joule heating
equations. We present optimal error estimates of
the semi-implicit Euler scheme in both the
$L^2$ norm and the $H^1$ norm without any time-step restriction.
Theoretical analysis is based on a new splitting of the error and
precise analysis of a corresponding time-discrete system.
The method used in this paper can be applied to more general nonlinear
parabolic systems and many
other linearized (semi)-implicit time discretizations for which
previous works often require certain restriction on the time-step size $\tau$.

\end{abstract}


\section{Introduction}
\setcounter{equation}{0} In the last several decades, numerous
effort has been devoted to the development of efficient numerical
schemes for nonlinear parabolic PDEs arising from a variety of
physical applications. A key issue to those schemes is the time-step
condition. Usually, fully implicit schemes are unconditionally
stable. However, at each time step, one has to solve a system of
nonlinear equations. An explicit scheme is much easy in computation.
But it suffers the severely restricted  time-step size for
convergence. A popular and widely-used approach is a linearized
(semi)-implicit scheme, such as linearized semi-implicit Euler
scheme. At each time step, the scheme only requires the solution of
a linear system. To study the error estimate of linearized
(semi)-implicit schemes, the boundedness of numerical solution (or
error function) in $L^{\infty}$ norm or a stronger norm is often
required. If a priori estimate for numerical solution in such a norm
cannot be obtained, one may employ the induction method with inverse
inequality to bound the numerical solution, such as
\begin{equation}
\| R_hu^n - U_h^n \|_{L^\infty} \le C h^{-d/2} \|R_h u^n - U_h^n \|_{L^2}
\le C h^{-d/2} (\tau^p + h^{r+1}) ,
\end{equation}
where $u^n$ and $U_h^n$ are the exact solution and numerical
solution, respectively, $R_h$ is some projection operator and $d$ is
the dimension. The above inequality, however, results in a time-step
restriction, particularly for problems in three dimensional space.
Such a technique has been widely used in error analysis for many
different nonlinear parabolic PDEs, $e.g.$, see \cite{AG,
He,KL,Liu1,Liu2} for Navier-Stokes equations, \cite{AL, EL, Zhao}
for nonlinear Joule heating problems, \cite{EW, SS,Wang} for porous
media flows, \cite{CL,EH, EM, WHS} for viscoelastic fluid flow,
\cite{MS,ZM} for KdV equations and \cite{DM,WML} for some other
equations. In all these works, error estimates were established
under certain time-step restrictions. We believe that these
time-step restrictions may not be necessary in most cases. In this
paper, we only focus our attention on a time-dependent and nonlinear
Joule heating system by a linearized semi-implicit scheme. However,
our approach is applicable for more general nonlinear parabolic PDEs
and many other time discretizations to obtain optimal error
estimates unconditionally.

The time-dependent nonlinear Joule heating system is defined by
\begin{align}
&\frac{\partial u}{\partial t}-\Delta
u=\sigma(u)|\nabla\phi|^2,
\label{e-heat-1}
\\[3pt]
&-\nabla\cdot(\sigma(u)\nabla\phi)=0,
\label{e-heat-2}
\end{align}
for $x\in\Omega$ and $t\in[0,T]$, where $\Omega$ is a bounded smooth
domain in $\R^d$, $d=2,3$. The initial and boundary conditions are given
by
\begin{align}
\label{BC}
\begin{array}{ll}
u(x,t)=0,\quad \phi(x,t)=g(x,t)~~
&\mbox{for}~~x\in\partial\Omega,~~t\in[0,T],\\[3pt]
u(x,0)=u_0(x)~~ &\mbox{for}~~x\in\Omega .
\end{array}
\end{align}

The nonlinear system above describes the model of electric heating
of a conducting body, where $u$ is the temperature, $\phi$ is the
electric potential, and $\sigma$ is the temperature-dependent
electric conductivity. Following the previous works \cite{EL,Zhao},
we assume that $\sigma\in W^{1,\infty}(\R)$ and
\begin{align}
\label{sigma}
\kappa\leq \sigma(s)\leq K,
\end{align}
for some positive constants $\kappa$ and $K$.

Theoretical analysis for the Joule heating system was done by
several authors \cite{AX,ALM,Cim, Xie, Yuan1,Yuan2, YL}. Among
these works, Yuan \cite{YL} proved existence and uniqueness of a
$C^\alpha$ solution in three-dimensional space. Based on this
result, further regularity can be derived with suitable
assumption on the initial and boundary conditions.
Numerical methods and
analysis for the Joule heating system can be found in \cite{AL, AY,
EL, Yue, Zhao, Zhou, ZW}. For the system in two-dimensional space,
optimal $L^2$ error estimate of a mixed finite element method with
the linearized semi-implicit Euler scheme was obtained in
\cite{Zhao} under a weak time-step condition. Error analysis for the
three-dimensional model was given in \cite{EL}, in which the
linearized semi-implicit Euler scheme with a linear Galerkin FEM was
used. An optimal $L^2$-error estimate was presented under the time
step restriction $\tau = O(h^{d/6})$. A more general time
discretization with higher-order finite element approximations was
studied in \cite{AL}. An optimal $L^2$-norm error estimate was given
under the conditions $\tau = O(h^{d/2p})$ and $r \ge 2$ where $p$
is the order of the discrete scheme in time direction and $r$ is the
degree of piecewise polynomial approximations used.
No optimal error estimates in $H^1$-norm have been obtained.

The main idea in this paper is a splitting of the numerical error
into the temporal direction and the spatial direction by introducing
a corresponding time-discrete parabolic system (or elliptic system).
Error bounds of the Galerkin finite element methods for the
time-discrete parabolic equations in certain norm is dependent only
upon the spatial mesh size $h$ and independent of the time-step size
$\tau$. If a suitable regularity of the solution of the
time-discrete equations can be proved, numerical solution in the
$L^{\infty}$ norm (or stronger norm) is bounded unconditionally
by the induction assumption together with the inverse inequality
\begin{equation}
\| R_hU^n - U_h^n \|_{L^\infty} \le C h^{-d/2} \| R_hU^n - U_h^n \|_{L^2}
\le C h^{r+1-d/2} ,
\end{equation}
where $U^n$ is the solution of the time-discrete equations.
With the boundedness,  optimal error estimates can
be established for the fully discrete scheme without any time-step
restriction. In this paper, we analyze the linearized
(semi-implicit) backward Euler scheme with the standard Galerkin
approximation in spatial directions for the nonlinear Joule heating
system \refe{e-heat-1}-\refe{BC}. With the splitting, we
present unconditionally optimal error estimates in
both the $L^2$ norm and the $H^1$ norm.

The rest of the paper is organized as follows. In Section 2, we
present the linearized semi-implicit Euler scheme with a linear
Galerkin finite element approximation in spatial direction and our
main results. After introducing the corresponding time-discrete
parabolic system, we provide in Section 3 a priori estimates and
optimal error estimates for the time-discrete solution, which imply
the suitable regularity of the time-discrete solution. With the
regularity obtained, we present optimal error estimates of the
Galerkin finite element solution in $L^2$-norm without any time-step
restriction, and the optimal error estimate in $H^1$ norm follows
immediately due to the nature of our approach. The concluding
remarks are presented in Section 4. Extension to $r$-order Galerkin
finite element approximation is straightforward with the
corresponding assumptions of regularity.

\section{Galerkin methods and main results}\label{dsfkjluweowefhiowheo}
\setcounter{equation}{0}
Let $\Omega$ be a bounded convex and smooth domain in $\R^d$ $(d=2,3)$. For any
integer $k\geq 0$ and $1\leq p<\infty$. Let
$W^k_p(\Omega)$ be the Sobolev space with the norm
$$
\|f\|_{W^k_p}=\biggl(\sum_{|\beta|\leq k}\int_\Omega|D^\beta f|^p\d
x\biggl)^\frac{1}{p},
$$
where
$$
D^\beta=\frac{\partial^{|\beta|}}{\partial
x_1^{\beta_1}\cdots\partial x_d^{\beta_d}}
$$
for the multi-index $\beta=(\beta_1,\cdots,\beta_d)$, $\beta_1\geq
0$, $\cdots$, $\beta_d\geq 0$, and $|\beta|=\beta_1+\cdots+\beta_d$.
For any integer $k\geq 0$ and $0<\alpha<1$, let
$C^{k+\alpha}(\overline\Omega)$ denote the usual H\"{o}lder space
with the norm
$$
\|f\|_{C^{k+\alpha}}=\sum_{|\beta|\leq k}\|D^\beta
f\|_{C(\overline\Omega)}+\sum_{|\beta|=
k}\sup_{x,y\in\Omega}\frac{|D^\beta f(x)-D^\beta
f(y)|}{|x-y|^\alpha}
$$
and let $C_0(\overline\Omega)$ be the space of continuous functions on $\overline\Omega$ vanishing on the boundary $\partial\Omega$.
For any Banach space $X$ and function $f:[0,T]\rightarrow X$,
we define the norm
$$
\|f\|_{L^p((0,T);X)} =\left\{
\begin{array}{ll}
\displaystyle\biggl(\int_0^T\|f(t)\|_X^pdt\biggl)^\frac{1}{p}, &
1\leq p<\infty
,\\[10pt]
\displaystyle{\rm ess\,sup}_{t\in(0,T)}\|f(t)\|_X, & p=\infty.
\end{array}
\right.
$$

With the boundary conditions (\ref{BC}), the weak
formulation of the system (\ref{e-heat-1})-(\ref{e-heat-2})
is defined by
\begin{eqnarray}
&& ( u_t, \, \xi_u) + ( \nabla u, \, \nabla \xi_u) = (\sigma(u)| \nabla
\phi |^2, \, \xi_u),
\label{weak-1}
\\
&& ( \sigma(u) \nabla \phi, \, \nabla \xi_{\phi} ) = 0
\label{weak-2}
\end{eqnarray}
for any $\xi_u,\xi_\phi\in H^1_0(\Omega)$ and a.e. $t\in(0,T)$.

Let $\pi_h$ be a regular division of $\Omega$ into triangles  in $\R^2$ or tetrahedras in $\R^3$,
i.e. $\Omega=\cup_{j}\Omega_j$, and denote by
$h=\max_{j}\{\mbox{diam}\,\Omega_j\}$ the mesh size. For a triangle $\Omega_j$ at the
boundary, we define $\widetilde\Omega_j$ as the triangle with one curved
side (or a tetrahedra with one curved face in $\R^3$) with the same
vertices as $\Omega_j$, and set $D_j=\widetilde \Omega_j\backslash
\Omega_j$. For an interior triangle, we set $\widetilde \Omega_j=\Omega_j$ and $D_j=\emptyset$. For a given division $\pi_h$, we define the finite element spaces \cite{VThomee}:
\begin{align*}
&V_{h}=\{v_h\in C(\overline\Omega): v_h|_{\Omega_j}\mbox{~is~linear~for~each~element~and~}v_h=0~\mbox{on}~ D_j\} ,\\
&S_{h}=\{v_h\in C(\overline\Omega): v_h|_{\widetilde \Omega_j}
\mbox{~is~linear~for~each~element}\} .
\end{align*}
It follows that $V_h$ is a subspace of $H^1_0(\Omega)$ and $S_h$ is
a subspace of $H^1(\Omega)$. For any function $v\in S_h$, we define $\Lambda_hv$ as the
function which satisfies $\Lambda_hv=0$ on $D_j$ and $\Lambda_hv=v$ on $T_j$. We define $\widetilde\Pi_h: C(\overline\Omega)\rightarrow
S_h$ to be the Lagrangian interpolation operator, i.e. $\widetilde\Pi_h v$ coincides with $v$ at each vertex of the triangular division of $\Omega$, and set
$\Pi_h=\Lambda_h\widetilde\Pi_h$. Clearly, $\Pi_h$ is a projection
operator from $C_0(\overline\Omega)$ onto $V_h$.

Let $\{ t_n \}_{n=0}^N$
be a partition in the time direction with $t_n = n \tau$, $T=N\tau$ and
$$
u^n = u(x,t_n), \quad \phi^n = \phi(x,t_n) \, .
$$
For any sequence of functions $\{ f^n \}_{n=0}^N$, we define
$$
D_t f^{n+1}=\frac{f^{n+1}-f^n}{\tau}\, .
$$

For simplicity, we assume that $g\in H^1(\Omega)$.
The fully discrete finite element scheme is to find
$U_h^n,\ \Phi_h^n
- g^n \in V_h$ for $n=0,1,\cdots,N$
such that for all $\xi_u,\ \xi_{\phi}\in V_h$
\begin{align}
&\big(D_t U_h^{n+1}, \, \xi_u \big)
+\big( \nabla U^{n+1}_h, \, \nabla \xi_u \big)
=\big(\sigma(U^n_h)|\nabla\Phi^n_h|^2, \, \xi_u \big),
\label{e-FEM-1}\\[3pt]
&\big(\sigma(U^{n}_h) \nabla\Phi^{n}_h, \, \nabla \xi_{\phi} \big)=0,
\label{e-FEM-2}
\end{align}
with the initial conditions $U_h^0 = I_h u^0$, where $I_h$ is the Lagrangian interpolation operator.

In the rest part of this paper, we always assume that the solution to the
initial/boundary value problem
(\ref{e-heat-1})-(\ref{BC}) exists and satisfies
\begin{align}
\label{StrongSOlEST}
&\|u\|_{L^\infty((0,T);H^2)}+\|u_t\|_{L^\infty((0,T);L^2)}+
\|u_{t}\|_{L^2((0,T);H^{2})} +\|u_{tt}\|_{L^2((0,T);L^2)} +\|u_0\|_{H^2} \nonumber\\
&
+\|\phi\|_{L^\infty((0,T);W^{2,12/5})}+\|\phi_t\|_{L^2((0,T);H^1)}+\|\nabla\phi\|_{L^\infty((0,T);C^\alpha)} \nonumber\\
&+\|g\|_{L^\infty((0,T);W^{2,12/5})}+\|g_t\|_{L^2((0,T);H^1)}+\|\nabla g\|_{L^\infty((0,T);C^\alpha)} \leq
C .
\end{align}
We denote by $C$ a generic positive
constant, which is independent of $n$, $h$ and $\tau$ and
$\epsilon$ a generic small positive constant.
We present our main results in the following theorem.

\begin{theorem}\label{ErrestFEMSol}
{\it
Suppose that the system \refe{e-heat-1}-\refe{e-heat-2}
with the initial and boundary conditions \refe{BC} has a unique
solution $(u, \phi)$ satisfying \refe{StrongSOlEST}. Then there
exist positive constants $\tau_0$ and $h_0$ such that when
$\tau<\tau_0$ and $h<h_0$, the finite element system
(\ref{e-FEM-1})-(\ref{e-FEM-2}) admits a unique solution $(U^n_h, \,
\Phi^n_h)$, $n=1,\cdots,N$, such that
\begin{align}
&\max_{1\leq n\leq N}\|U^n_h - u^n\|_{L^2} +\max_{1\leq n\leq
N}\|\Phi^n_h -\phi^n\|_{L^2} \leq C(\tau+h^2) , \label{optimalL2est} \\
&\max_{1\leq n\leq N}\|U^n_h - u^n\|_{H^1} +\max_{1\leq n\leq
N}\|\Phi^n_h - \phi^n\|_{H^1} \leq C(\tau+h). \label{optimalH1est}
\end{align}

}
\end{theorem}
\medskip

For $U^0=u_0$ and $\Phi^0$, we define
$U^n$ and $\Phi^n$ to be the solution of the following
discrete parabolic system (or elliptic system)
\begin{align}
&D_t U^{n+1} - \Delta U^{n+1}=\sigma(U^n)|\nabla \Phi^n|^2,\quad 0\leq n\leq N-1,
\label{TDEqjoulheat}\\[3pt]
&-\nabla\cdot(\sigma(U^n)\nabla \Phi^n)=0,\quad\quad\quad\quad\quad~~ 0\leq n\leq N,
\label{TDEqjoulheat23}
\end{align}
with the boundary conditions
\begin{align}\label{TDBDCsjoulheat}
\begin{array}{ll}
U^{n+1}(x)=0,\quad \Phi^n(x)=g(x,t_n)~~
&\mbox{for}~~x\in\partial\Omega .
\end{array}
\end{align}

We will present the proof of Theorem \ref{ErrestFEMSol}
in the next two sections.
The key to our proof is the following error splitting
\begin{align}
&\|U^n_h- u^n\| \le \| e^n \| + \| e_h^n \| + \| U^n - R_h U^n \|,
\nn \\
&\|\Phi_h^n- \phi^n\| \le \| \eta^n \|+ \| \eta_h^n \|
+ \| \Phi^n - P_h^n\Phi^n \|
\nn
\end{align}
for any norm $\|\cdot\|$, where
\begin{eqnarray}
&& e^n = U^n - u^n, \quad e_h^n = U_h^n - R_h U^n, \quad
\nn \\
&& \eta^n = \Phi^n - \phi^n, \quad \eta_h^n = \Phi_h^n - P_h^n \Phi^n \, ,
\nn
\end{eqnarray}
with
$R_h:H^1_0(\Omega)\rightarrow V_h$ being the Riesz projection
operator, i.e.
\begin{align*}
\!\!\!\!\!\!\!\!\!\!\!\!\!\big(\nabla (v-R_hv),\nabla
w\big)=0,\quad\mbox{for~all}~~v\in H^1_0(\Omega)~~\mbox{and}~~w\in
V_h.
\end{align*}  and $P_h^n\Phi^n=g(\cdot,t_n)+\Pi_h(\Phi^n-g(\cdot,t_n))$ for $n=0,1,2,\cdots,N$.

With the definition of the operator $\Pi_h$, $P_h^n$ and $R_h$,  the following estimates hold \cite{RS}:
for any $2\leq p<\infty$, there exists a positive constant $C$ (independent of the function $v$) such that
\begin{align}
&\|v-\Pi_hv\|_{L^p}+h\|v-\Pi_hv\|_{W^{1,p}}\leq
Ch^2\|v\|_{W^{2,p}},\label{ip-1}
\\[5pt]
&\|\Phi^n-P_h^n\Phi^n\|_{L^p}+h\|\Phi^n-P_h^n\Phi^n\|_{W^{1,p}} \leq
Ch^2\|\Phi^n-g^n\|_{W^{2,p}}, \\[5pt]
&\|v-R_hv\|_{L^p} +h\|v-R_hv\|_{W^{1,p}}\leq Ch^2 \|v\|_{W^{2,p}},
\label{ip-2}
\end{align}
for $v\in W^{2,p}(\Omega)\cap H^1_0(\Omega)$ .

\section{Error estimates}
\setcounter{equation}{0} We analyze the error function $(e^n,
\eta^n)$ from the linearized semi-implicit Euler scheme
(time-discrete system) and the errors function $(e^n_h, \eta^n_h)$
of the Galerkin finite element method for the time-discrete system
in the following two subsections, respectively.

\subsection{The time-discrete solution}
In this subsection, we prove the existence and uniqueness of the
time-discrete system (\ref{TDEqjoulheat})-(\ref{TDBDCsjoulheat})
and establish the error bounds for
$(e^n, \eta^n)$.

\begin{theorem}\label{ErrestDisSol}
{\it Suppose that the system
\refe{e-heat-1}-\refe{BC} has a unique solution $(u, \phi)$ satisfying
\refe{StrongSOlEST}.
Then there exists a positive constant $\tau_0$
such that when $\tau<\tau_0$, the time-discrete system
(\ref{TDEqjoulheat})-(\ref{TDBDCsjoulheat}) admits a unique solution
$(U^n, \Phi^n)$ such that
\begin{align}\label{ErrestDisSol222}
&\max_{1\leq n\leq N}\| U^n\|_{H^2}
+\max_{1\leq n\leq N}\|D_t U^n\|_{L^2}
+ \left ( \sum_{n=1}^{N} \tau \| D_t U^n\|_{H^2}^2 \right )^{1/2} \leq C,\\
&\max_{1\leq n\leq N}\| \Phi^n\|_{W^{2,12/5}} +\max_{1\leq n\leq
N}\|\nabla\Phi^n\|_{L^\infty}\leq C
\label{ErrestDisSol2223}
\end{align}
and
\begin{align}\label{TDErrEstLemm}
\begin{array}{ll}
&\displaystyle\max_{1\leq n\leq N}\| e^{n}\|_{H^1} +\max_{1\leq
n\leq N}\| \eta^{n}\|_{H^1} \leq C \tau .
\end{array}
\end{align}
}
\end{theorem}

\noindent{\it Proof}~~~ We rewrite the system \refe{e-heat-1}-\refe{e-heat-2}
by
\begin{align}
\label{Eqjoulheat2} &D_t u^{n+1}-\Delta u^{n+1}=\sigma( u^n)|\nabla
\phi^n|^2
+R_1^{n+1},\\[3pt]
&-\nabla\cdot(\sigma( u^{n})\nabla\phi^{n})=0,
\label{Eqjoulheat22}
\end{align}
where $R_1^{n+1}$ is the truncation errors due to the time
discretization, i.e.
\begin{align*}
R_1^{n+1}=D_tu^{n+1}-\frac{\partial u}{\partial
t}\Big|_{t=t_{n+1}}+(\sigma(u^{n+1})-\sigma(u^n))|\nabla\phi^{n+1}|^2\\
+\sigma(u^n)\nabla(\phi^{n+1}+\phi^n)\cdot\nabla(\phi^{n+1}-\phi^n)
.
\end{align*}
With the regularity given in (\ref{StrongSOlEST}), we have
\begin{align}\label{truncerr}
\begin{array}{ll}
&
\|R_1^{n+1}\|_{L^2}\leq C , \quad
\sum_{n=0}^{N-1}\|R_1^{n+1}\|_{L^2}^2 \tau \leq C \tau^2 .
\end{array}
\end{align}
Subtracting the equations \refe{Eqjoulheat2}-\refe{Eqjoulheat22}
from the equations (\ref{TDEqjoulheat})-(\ref{TDEqjoulheat23}),
respectively, we obtain
\begin{align}
&D_t e^{n+1}-\Delta e^{n+1}=(\sigma( U^n) - \sigma(u^n))
|\nabla \phi^n|^2 \nonumber\\
&~~~~~~~~~~~~+\sigma(U^n)(\nabla\phi^n +\nabla \Phi^n)\cdot\nabla
\eta^n+R_1^{n+1},
\label{TDErrEq1}\\[5pt]
&-\nabla\cdot(\sigma(U^{n}) \nabla \eta^{n})
=\nabla\cdot[(\sigma( u^{n})-\sigma(U^{n}))\nabla\phi^{n}] \, .
\label{TDErrEq2}
\end{align}
An alternative to the last equation is
\begin{align}
&-\nabla\cdot(\sigma(u^{n})\nabla \eta^{n})
=\nabla\cdot[(\sigma(u^{n})-\sigma(U^{n}))(\nabla\phi^{n}
+\nabla\eta^{n})]
 .
\label{TDErrEq3}
\end{align}
Multiplying the equation (\ref{TDErrEq2}) by $\eta^{n+1}$ and integrating
the result over $\Omega$, we have
$$
\|\nabla \eta^{n}\|_{L^2}^2\leq C\|e^{n} \|_{L^2}
\| \nabla \eta^{n}\|_{L^2}
$$
which leads to
\begin{equation}
\|\nabla \eta^{n}\|_{L^2}\leq C\|e^{n} \|_{L^2} .
\label{TDErrEq2est}
\end{equation}

Similarly,  multiplying (\ref{TDErrEq1}) by $e^{n+1}$ and integrating it
over $\Omega$ gives
\begin{align*}
&D_t\biggl(\frac{1}{2}\|e^{n+1}\|_{L^2}^2\biggl)+\|\nabla e^{n+1}\|_{L^2}^2
 \\
&\leq C\|e^n\|_{L^2} \|e^{n+1} \|_{L^2} \|\nabla \phi^n
\|_{L^{\infty}} + \big(\sigma(U^n)(\nabla\phi^n +\nabla \Phi^n)
e^{n+1}, \, \nabla \eta^n\big)
\\
& \, \, + \|R_1^{n+1}\|_{L^2} \| e^{n+1} \|_{L^2} \, .
\end{align*}
By \refe{TDEqjoulheat23} and using integrating by part,
\begin{align*}
& |\big(\sigma(U^n)(\nabla\phi^n +\nabla \Phi^n)e^{n+1}, \,
\nabla \eta^n\big)|\\
& \qquad \le |\big(\sigma(U^n)e^{n+1}\nabla\phi^n, \, \nabla
\eta^n\big)|
\\
& \qquad \quad + | \big(e^{n+1} \nabla \cdot (\sigma(U^n) \nabla
\Phi^n) + \sigma(U^n) \nabla \Phi^n \cdot \nabla e^{n+1}, \, \eta^n
\big)|
\\
&\qquad \le C (\| e^{n+1} \|_{L^2} \|\nabla\eta^n \|_{L^2} +
\|\nabla \phi^n\|_{L^\infty}\|\nabla
e^{n+1} \|_{L^2}\|\eta^n\|_{L^2}+\|\nabla
\eta^n\|_{L^2}\|\nabla e^{n+1} \|_{L^2}\|\eta^n\|_{L^\infty} )
\, .
\end{align*}
Applying the maximum principle to the elliptic equation
(\ref{TDEqjoulheat23}) shows that $\| \Phi^{n} \|_{L^{\infty}} \leq C$
and therefore,
$$
\| \eta^n \|_{L^{\infty}} \leq C ,
$$
for $n=0,1,2,\cdots$. It follows that
\begin{align*}
&D_t\biggl(\frac{1}{2}\|e^{n+1}\|_{L^2}^2\biggl)+\frac{1}{2}\|\nabla
e^{n+1}\|_{L^2}^2 \\
&\leq C\|e^n\|_{L^2}^2+C\| e^{n+1}\|_{L^2}^2
+C\|\eta^n\|_{H^1}^2+C\|R_1^{n+1}\|_{L^2}^2\\
&\leq C\| e^n\|_{L^2}^2+C\|e^{n+1}\|_{L^2}^2 +C\|R_1^{n+1}\|_{L^2}^2 ,
\end{align*}
where we have used (\ref{TDErrEq2est}) in the last step.
By applying Gronwall's inequality, combined with
(\ref{truncerr}), we derive that there exists a small positive constant $\tau_0$ such that when $\tau<\tau_0$,
\begin{align}\label{nksfdljiweo}
&\max_{1\leq n\leq N}\| e^{n}\|_{L^2}^2+\max_{1\leq n\leq
N}\| \eta^{n}\|_{H^1}^2+\sum_{n=1}^N\| e^{n}\|_{H^1}^2 \tau \leq
C \tau^2.
\end{align}
In particular, the above estimate implies that
\begin{align}\label{estibaruh1}
&\|U^n\|_{H^1}^2\leq C
\end{align}
and
$$
\|D_tU^{n+1}\|_{L^2}\leq
\|D_tu^{n+1}\|_{L^2}+\|D_te^{n+1}\|_{L^2}\leq C .
$$
With the above inequalities, we derive from (\ref{TDEqjoulheat})
that
\begin{align}\label{H2estintermsPhi}
\|U^{n+1}\|_{H^2} & \leq C+C\|\nabla\Phi^n\|_{L^4}^2 .
\end{align}
Since $H^2(\Omega)\hookrightarrow C^\alpha(\overline\Omega)$ in
$\R^d$ with $d=2,3$, $\| e^n \|_{C^{\alpha}} \le C$. By applying the
$W^{1,4}$ estimate \cite{ByunWang,Simader} to (\ref{TDErrEq3}), we
get
\begin{align*}
\|\nabla\eta^{n}\|_{L^4} &\leq \| (\sigma(u^{n}) - \sigma(U^{n}))
\nabla \Phi^{n} \|_{L^4}
\\
& \leq C_0\| e^n\|_{L^\infty}( \|\nabla \phi^{n} \|_{L^4} + \|\nabla
\eta^{n}\|_{L^4})
\end{align*}
where $C_0$ is some positive constant. By assuming that $C_0\|
e^n\|_{L^\infty}<1/2$, we derive that
\begin{align*}
\|\nabla\eta^{n}\|_{L^4} \leq C
\end{align*}
and (\ref{H2estintermsPhi}) implies that
\begin{align}\label{H2estabs}
\| e^{n+1} \|_{H^2} \le \| u^{n+1} \|_{H^2} + \|U^{n+1}\|_{H^2}\leq C
\end{align}
and
$$
\| e^{n+1}\|_{L^\infty} \leq \| e^{n+1}\|_{H^1}^{1/2}\|
e^{n+1}\|_{H^2}^{1/2} \leq C\tau^{1/4} .
$$
From the above derivation, one can see that there exists $\tau_0>0$
such that if $\tau<\tau_0$, then $C_0\| e^n\|_{L^\infty}<1/2$
implies $C_0\| e^{n+1}\|_{L^\infty}<1/2$ as well as
(\ref{H2estabs}).
In addition, we see that $\|\nabla \Phi^n \|_{L^4}\leq C$ and
therefore,
\begin{align}
\max_{1\leq n\leq N}\| U^{n}\|_{C^\alpha}\leq C .
\end{align}

By applying Schauder's estimates (\cite{ChenYZ}, page 74) to
(\ref{TDEqjoulheat23}), we derive that
\begin{align}\label{estphiCalpha}
\max_{1\leq n\leq N}\|\nabla\Phi^{n}\|_{C^\alpha}\leq C ,
\end{align}
which together with (\ref{estibaruh1}) and (\ref{TDEqjoulheat23})
implies that
\begin{align}\label{estphiH2}
\max_{1\leq n\leq N}\|\Phi^{n}\|_{W^{2,12/5}}\leq C .
\end{align}

Multiplying  (\ref{TDErrEq1}) by $-\Delta e^{n+1}$ and summing up
the equations for $n=0,1,\cdots,N-1$, we obtain
\begin{align*}
&\max_{1\leq n\leq N}\| e^{n}\|_{H^1}^2+\sum_{n=1}^N\tau\|\Delta e^n\|_{L^2}^2\\
&\leq \sum_{n=0}^{N-1}\tau\Big(\| (\sigma( U^n) - \sigma(u^n)) |\nabla
\phi^n|^2 \|_{L^2}^2
+ \| \sigma(U^n)(\nabla\phi^n +\nabla \Phi^n)\cdot\nabla \eta^n
\|_{L^2}^2 + \| R_1^{n+1}\|_{L^2}^2\Big)
\leq C\tau^2 .
\end{align*}
It follows that
\begin{align*}
&\max_{1\leq n\leq N}\| e^{n}\|_{H^1}\leq C\tau,
\end{align*}
and
\begin{align*}
&\sum_{n=1}^N\tau\|\Delta D_te^n\|_{L^2}^2
\leq C\tau^{-2}\sum_{n=1}^N\tau\|\Delta e^n\|_{L^2}^2\leq C .
\end{align*}
By the theory of elliptic equations \cite{ChenYZ, Evans},
$\|D_te^n\|_{H^2}\leq C\|\Delta D_te^n\|_{L^2}$ for $n=1,\cdots,N$, and so
\begin{align} \label{TDEqjoulheat33}
&\sum_{n=1}^N\tau\|D_te^n\|_{H^2}^2\leq  C .
\end{align}

The proof of Theorem \ref{ErrestDisSol} is complete. \endproof

\subsection{The fully-discrete finite element solution}
Here we study the error $(e^n_h, \eta^n_h)$
of the Galerkin finite element method for the time-discrete system
(\ref{TDEqjoulheat})-(\ref{TDBDCsjoulheat}).

\begin{theorem}\label{full-error}
Suppose that the system
\refe{e-heat-1}-\refe{BC} has a solution $(u, \phi)$ satisfying
\refe{StrongSOlEST}.
Then there exist positive constants $h_0$ and $\tau_0$
such that when $h<h_0$ and $\tau<\tau_0$, the fully-discrete system
\refe{e-FEM-1}-\refe{e-FEM-2} admits a unique solution
$(U_h^n, \Phi^n_h)$ such that
\begin{align}
& \| e^n_h\|_{L^2}
+ \| \eta^n_h\|_{L^2} \leq C h^2  , \label{full-error-2}\\
& \| e^n_h\|_{H^1}
+\|\eta_h^n \|_{H^1} \le Ch .
\label{full-error-1}
\end{align}
\end{theorem}

Note that the condition of $\tau<\tau_0$ is to ensure that Theorem
\ref{ErrestDisSol} holds. For the given $U^n_h$ , the error
estimate for the equation \refe{e-FEM-2} is given in the following Lemma.
\begin{lemma}\label{lemma4.3}
{\it Suppose that the system
\refe{e-heat-1}-\refe{BC} has a unique solution $(u, \phi)$ satisfying
\refe{StrongSOlEST}. Then
\begin{align*}
&\|\nabla(\Phi^{n}_h - \Phi^{n})\|_{L^2}
\leq C\big(h+\|e^{n}_h\|_{L^2}\big),\\
&\|\Phi^{n}_h - \Phi^{n}\|_{L^2}\leq
C\big(h^2+\|e^{n}_h\|_{L^2}+h^{-d/6}\|e^{n}_h\|_{L^2}^2\big)
,
\end{align*}
where $(U_h^n, \, \Phi_h^n)$ and
$(U^n, \, \Phi^n)$ are the solution of
the finite element system \refe{e-FEM-1}-\refe{e-FEM-2} and
the time-discrete system
(\ref{TDEqjoulheat})-(\ref{TDBDCsjoulheat}), respectively.
}
\end{lemma}

\noindent{\bf Remark 3.1}~\, The proof of the above lemma is similar
as that of Lemma 3.2 in \cite{EL}, in which the factor $h^{-d/6}$
appears when $\|e^n_h \|_{L^3}$ reduces to $\|e^n_h \|_{L^2}$ via
the inverse inequality. More important is that in \cite{EL}, $e_h^n$
is the difference between the exact solution of the system
\refe{e-heat-1}-\refe{e-heat-2} and the fully discrete finite
element solution. The restriction for the time-step size, $\tau \le
k_0 h^{d/6}$, was required when the preliminary error bound
$\|e_h^n\|_{L^2} \le C(\tau + h^2)$ was used by induction in the
second inequality of Lemma \ref{lemma4.3}. However, in our approach,
$e_h^n$ is the difference between the solution of the time-discrete
system (\ref{TDEqjoulheat})-(\ref{TDBDCsjoulheat}) and the fully
discrete finite element solution. Thus, the induction assumption
shows that $\|e_h^n\|_{L^2} \le C h^2$ and then, we can prove the
optimal error bound of the scheme unconditionally. \vskip0.1in

\noindent{\it Proof of Theorem \ref{full-error}}~~
At each time step of the scheme,
one only needs to solve two uncoupled linear discrete systems.
Due to the assumption \refe{sigma}, it is easy to see that coefficient matrices in both systems are symmetric
and positive definite.
Existence and uniqueness of the Galerkin finite element solution follows
immediately. It is seen that the inequality (\ref{full-error-1}) follows from (\ref{full-error-2}) via the inverse inequality. Therefore, it suffices to prove (\ref{full-error-2}).

The weak formulation of the time-discrete system
(\ref{TDEqjoulheat})-(\ref{TDBDCsjoulheat}) is
\begin{align}
&\big(D_t U^{n+1}, \, \xi_u \big) +\big( \nabla U^{n+1}, \, \nabla
\xi_u \big) =\big(\sigma(U^n)|\nabla\Phi^n|^2, \, \xi_u \big),
\label{d-FEM-1}\\[3pt]
&\big(\sigma(U^{n}) \nabla\Phi^{n}, \, \nabla \xi_{\phi}
\big)=0, \label{d-FEM-2}
\end{align}
for any $\xi_u,\ \xi_{\phi}\in V_h$. From the above equations and
the finite element system \refe{e-FEM-1}-\refe{e-FEM-2}, we find
that the error function $(e_h^n, \eta_h^n)$ satisfies
\begin{align}
&\big(D_te^{n+1}_h,\, \xi_u \big)+\big(\nabla e^{n+1}_h,\,
\nabla \xi_u \big)
\nn\\
&=\big(D_t(U^{n+1}-R_h U^{n+1}), \, \xi_u \big) +\big((\sigma(U^n_h)
- \sigma(U^n))|\nabla\Phi^n|^2, \, \xi_u \big)
\nn\\
&~~~~~+2\big((\sigma(U^n_h)-\sigma(U^n))
\nabla\Phi^n\cdot\nabla(\Phi^n_h-\Phi^n), \, \xi_u \big)
\nn\\
&~~~~~+\big(\sigma(U^n_h)|\nabla(\Phi^n_h-\Phi^n)|^2, \, \xi_u \big)
\nn\\
&~~~~~+2\big(\sigma(U^n)\nabla\Phi^n\cdot
\nabla(\Phi^n_h-\Phi^n), \, \xi_{u} \big)\nonumber\\
&:=(\bar R_1^{n+1}, \, \xi_u )+(\bar R_2^{n+1}, \, \xi_u)+(\bar
R_3^{n+1},\, \xi_u )+(\bar R_4^{n+1}, \, \xi_u)+(\bar
R_5^{n+1}, \, \xi_u),
\label{FEMErrEq1}
\end{align}
and
\begin{align}
& \big(\sigma(U^{n}_h)\nabla\eta^{n}_h, \, \nabla \xi_{\phi}
\big) =\big((\sigma(U^{n})-\sigma(U^{n}_h))\nabla \Phi^{n}, \,
\nabla \xi_{\phi} \big)
\nn\\
&~~~~~~~~~~~~~~~~~~~~~~~~~~~~~
+\big(\sigma(U^{n}_h)\nabla(\Phi^{n}-P^{n}_h\Phi^{n}), \, \nabla
\xi_{\phi} \big) \label{FEMErrEq2}
\end{align}
for all $\xi_u, \xi_{\phi}\in V_h$.


Since $\eta^{n}_h=0$ on $\partial\Omega$,  we can take
$\xi_{\phi}=\eta^{n}_h$ in (\ref{FEMErrEq2}) to get
\begin{align}\label{estdetal21}
&\|\nabla\eta^{n}_h\|_{L^2}\leq C\| e^{n}_h\|_{L^2}+Ch,
\end{align}
where we have noted the fact that
$\|\nabla(\phi^{n}-P^{n}_h\phi^{n})\|_{L^2}\leq Ch$. With the
above inequality, from Lemma \ref{lemma4.3} we derive that
\begin{align}\label{estdetal2133}
\|\eta^n_h\|_{L^2}\leq Ch^2+C\| e^n_h\|_{L^2} +h^{-d/6}\|e^{n}_h\|_{L^2}^2.
\end{align}

Taking $\xi_u=e^{n+1}_h$ in (\ref{FEMErrEq1}), the right-hand side
is estimated by
\begin{align}
(\bar R_1^{n+1}, \, e^{n+1}_h) & \leq \epsilon\|e^{n+1}_h\|_{H^1}^2
+C\epsilon^{-1}\|D_t U^{n+1}-R_h D_t U^{n+1}\|_{L^2}^2 \nn\\
&\leq \epsilon\|e^{n+1}_h\|_{H^1}^2+C\epsilon^{-1}\|D_t
U^{n+1}\|_{H^2}^2h^4,
 \label{R1herrestimate}\\
(\bar R_2^{n+1}, \, e^{n+1}_h) & \leq C\|e^{n+1}_h\|_{L^2} (
\|e^{n}_h \|_{L^2} + \| U^n - R_h U^n \|_{L^2} )
 \nn\\
& \leq \epsilon \|e^{n+1}_h\|_{L^2}^2 + C\epsilon^{-1}(\|e^{n}_h \|_{L^2}^2 + h^4),
 \label{R2herrestimate}\\
(\bar R_3^{n+1}, \, e^{n+1}_h) & \leq C \|e^{n+1}_h\|_{L^6}(
\|e^n_h\|_{L^2} + \| U^n - R_h U^n \|_{L^2} ) (
\|\nabla\eta^{n+1}_h\|_{L^3} + \| \Phi^n - P^n_h\Phi^n \|_{L^3})
 \nn\\
& \leq C\|e^{n+1}_h\|_{H^1}(\|e^n_h\|_{L^2}+Ch^2)
(\|\nabla\eta^n_h\|_{L^3}+Ch)
 \nn\\
&\leq \epsilon\|e^{n+1}_h\|_{H^1}^2+
C\epsilon^{-1}(\|e^n_h\|_{L^2}+Ch^2)^2
(Ch^{-d/6}\|\nabla\eta^n_h\|_{L^2}+Ch)^2
 \nn\\
&\leq \epsilon\|e^{n+1}_h\|_{H^1}^2+
C\epsilon^{-1}(\|e^n_h\|_{L^2}^2+h^4)(Ch^{-d/6}\|e^n_h\|_{L^2}+Ch^{1-d/6})^2,
 \label{R3herrestimate}\\
(\bar R_5^{n+1}, \, e^{n+1}_h)
&=-2\big(\sigma(U^n)\nabla\Phi^n\cdot
(\Phi^n_h-\Phi^n), \, \nabla e^{n+1}_h\big)\\
& \leq
C\|\Phi^n_h-\Phi^n\|_{L^2}\|\nabla e^{n+1}_h\|_{L^2}
 \nn\\
& \leq \epsilon\| \nabla e^{n+1}_h\|_{L^2}^2
+C\epsilon^{-1}(\|e^n_h\|_{L^2}^2+h^{-d/3}\|e^n_h\|_{L^2}^4+h^4 ) \label{R5herrestimate}
\end{align}
and
\begin{align}
(\bar R_4^{n+1},\, e^{n+1}_h) & \leq C \|e^{n+1}_h\|_{L^\infty}
( \|\nabla\eta^n_h\|_{L^2}^2
+ \| \nabla (\Phi^n - P^n_h\Phi^n) \|_{L^2}^2 )
 \nn\\
& \leq  C h^{-1/2} \|e^{n+1}_h\|_{H^1} (
\|e^n_h\|_{L^2}^2 + h^2)
\nn \\
& \leq \epsilon \|e^{n+1}_h\|_{H^1}^2 + C h^{-1} \| e_h^n \|_{L^2}^4 +
C h^3 \, .
\label{R4est}
\end{align}
With the above estimates, (\ref{FEMErrEq1}) reduces to
\begin{align}\label{fsdklfjweiof}
D_t \left ( \|e^{n+1}_h\|_{L^2}^2 \right )
+ \| \nabla e_h^{n+1} \|_{L^2}^2 \leq & C \left ( \|e^{n}_h\|_{L^2}^2
+ h^{-1} \|e^{n}_h\|_{L^2}^4 \right ) \nn\\
& + C h^3+C\epsilon^{-1}\|D_t
U^{n+1}\|_{H^2}^2h^4 ,
\,
\end{align}
which holds for $0\leq n\leq N-1$.

Now we prove that
\begin{align}\label{dfskweuio}
\|e^n_h\|_{L^2}\leq h^{1/2}~~\mbox{for}~~0\leq n\leq N
\end{align}
by using mathematical induction. Clearly, this inequality holds for $n=0$. If we assume that this inequality holds for $0\leq n\leq k$, then the inequality (\ref{fsdklfjweiof}) reduces to
\begin{align}\label{fsdklfjweiof1}
D_t \left ( \|e^{n+1}_h\|_{L^2}^2 \right )
+ \| \nabla e_h^{n+1} \|_{L^2}^2 \leq & C  \|e^{n}_h\|_{L^2}^2+ C h^3+C\epsilon^{-1}\|D_t
U^{n+1}\|_{H^2}^2h^4
\,
\end{align}
for $0\leq n\leq k$. By applying Gronwall's inequality, we derive that
\begin{align}\label{fsdklfjweiof2}
\|e^{k+1}_h\|_{L^2}^2\leq C_{1}h^3\leq h^{1/2} ~~\mbox{if}~~h<1/C_{1}^{2/5} .
\end{align}
This completes the induction.

With (\ref{dfskweuio}), we can apply Gronwall's inequality to
(\ref{fsdklfjweiof}) and get
\begin{align}\label{fsdklfjweiof2}
\max_{1\leq n\leq N}\|e^n_h\|_{L^2}\leq Ch^{3/2}  .
\end{align}

Since $\eta^{n+1}_h \in H^1_0(\Omega)$,
from the estimates (\ref{estdetal21})-(\ref{estdetal2133})  we see that
\begin{align}
\max_{1\leq n\leq N}\|\eta^n_h \|_{H^1}\leq Ch .
\end{align}
which implies that $\|\nabla \eta_h^n\|_{L^6}\leq Ch^{-d/3}\|\nabla \eta_h^n\|_{L^2}\leq C$ and so $\|\nabla \Phi_h^n\|_{L^6}\leq C$.

Finally, we rewrite (\ref{FEMErrEq2}) as
\begin{align}
& \big(\sigma(U^{n})\nabla\eta^{n}_h, \, \nabla \xi_{\phi}
\big)\nn\\
& =\big((\sigma(U^{n})-\sigma(U^{n}_h))\nabla\Phi^{n}_h, \, \nabla \xi_{\phi}
\big)
+\big(\sigma(U^{n})\nabla(\Phi^{n}-P^{n}_h\Phi^{n}), \, \nabla
\xi_{\phi} \big) ,\quad\forall~\xi_\phi\in V_h ,
\end{align}
and apply
the $W^{1,p}$ estimate \cite{RS} to the above equation. Then we get
\begin{align}
\|\nabla\eta^{n}_h\|_{L^{12/5}}
&\leq C\|(\sigma(U^{n})-\sigma(U^{n}_h))\nabla\Phi^{n}_h\|_{L^{12/5}}
+C\|\nabla(\Phi^{n}-P^{n}_h\Phi^{n})\|_{L^{12/5}} \nn\\
&\leq C
\|U^{n}_h-U^{n}\|_{L^4}\|\nabla \Phi_h^n\|_{L^6}
+C\|\nabla(\Phi^{n}-P^{n}_h\Phi^{n})\|_{L^{12/5}}\nn\\
&\leq
Ch^{-d/4}(\|e^{n}_h\|_{L^2}+h^{2}\|U^{n}\|_{H^2})+Ch\|\Phi^{n}\|_{W^{2,12/5}}\nn\\
&\leq C\|e^{n}_h\|_{L^2}^{1/2}+Ch .
\end{align}
Therefore, we obtain a refined estimate:
\begin{align}\label{R4herrestimate}
(\bar R_4^{n+1},\, e^{n+1}_h) & \leq C \|e^{n+1}_h\|_{L^6}
( \|\nabla\eta^n_h\|_{L^{12/5}}^2
+ \| \nabla (\Phi^n - P^n_h \Phi^n) \|_{L^{12/5}}^2 )
 \nn\\
& \leq \epsilon \|e^{n+1}_h\|_{L^6}^2 + C\epsilon^{-1}
\|\nabla\eta^n_h\|_{L^{12/5}}^4 +
C\epsilon^{-1} \|\nabla(\Phi^n-P^n_h\Phi^n)\|_{L^{12/5}}^4
\nn \\
& \leq \epsilon \|e^{n+1}_h\|_{H^1}^2 +C\epsilon^{-1}
\|e^{n}_h\|_{L^2}^2 + C\epsilon^{-1} h^4 \, .
\end{align}
With the estimates (\ref{R1herrestimate})-(\ref{R5herrestimate}) and (\ref{R4herrestimate}), the equation (\ref{FEMErrEq1}) reduces to
\begin{align*}
D_t \left ( \|e^{n+1}_h\|_{L^2}^2 \right )
+ \| \nabla e_h^{n+1} \|_{L^2}^2 \leq  C  \|e^{n}_h\|_{L^2}^2+ C h^4+C\|D_t
U^{n+1}\|_{H^2}^2h^4 .
\,
\end{align*}
By applying Gronwall's inequality, we get
\begin{align}\label{dsfjkwojefiwo}
\max_{1\leq n\leq N}\|e^{n}_h\|_{L^2}^2    \leq Ch^4 .
\end{align}
The $L^2$ error estimate of $\eta_h^n$ follows from (\ref{estdetal2133}) and (\ref{dsfjkwojefiwo}).
The proof of Theorem \ref{full-error} is complete. \endproof
\vskip0.1in

Theorem \ref{ErrestFEMSol} follows immediately from Theorem \ref{full-error}
and Theorem \ref{ErrestDisSol}. \endproof
\vskip0.1in

\section{Conclusions}
\setcounter{equation}{0} We have presented an approach to obtain
optimal error estimates and unconditional stability of linearized
(semi) implicit schemes with a Galerkin finite element method for
the three-dimensional nonlinear Joule heating equations.  The
analysis is based on a new splitting of the error into the time
direction and the spatial direction, by which the numerical solution
(or its error) in a strong norm can be bounded by induction
assumption and the inverse inequalities without any restrictions on
the time-step size. In most existing approaches, a time-step
condition has to be enforced to bound the numerical solution in a
stronger norm. Clearly, our analysis can be extended to many other
nonlinear parabolic systems, while we only focus on the electric
heating model in the present paper.

In this paper, we only considered a linear Galerkin finite element
approximation. The extension to high-order Galerkin finite element
methods can be done similarly. For simplicity, we have assumed that
the function $g$ is defined in the domain $\Omega$ instead of on the
boundary $\partial\Omega$. If the function $g$ is defined only on
the boundary $\partial\Omega$, a similar analysis can be given by
taking the boundary terms into consideration, see \cite{EL} for
reference. Optimal error estimates still can be proved without any
condition on the time-step size.



\end{document}